\documentclass[reqno]{amsart}

\usepackage{amssymb}
\usepackage{amsfonts}

\newtheorem{theorem}{Theorem}
\theoremstyle{plain}

\newtheorem{definition}{Definition}

\newtheorem{lemma}{Lemma}

\newtheorem{remark}{Remark}

\numberwithin{equation}{section}
\input{tcilatex}

\begin{document}
\title[Kato Estimates for Harmonic Oscillator]{Classical Proofs Of Kato Type
Smoothing Estimates for The Schr\"{o}dinger Equation with Quadratic
Potential in $\mathbb{R}^{n+1}$ with Application}
\author{Xuwen Chen}
\address{Department of Mathematics, University of Maryland, College Park, MD
20742}
\email{chenxuwen@math.umd.edu}
\thanks{}
\date{07/01/2010}

\begin{abstract}
In this paper, we consider the Schr\"{o}dinger equation with quadratic
potential
\begin{equation*}
i\frac{\partial }{\partial t}u=-\triangle u+\left\vert x\right\vert ^{2}u
\text{ }in\text{ }\mathbb{R}^{n+1},\text{ }u(x,0)=f(x)\in L^{2}(\mathbb{R}
^{n}).
\end{equation*}
Using Hermite functions and some other classical tools, we give an
elementary proof of the Kato type smoothing estimate: for $i\neq j\neq k,$ $
\delta \in \lbrack 0,1],$ and $n\geqslant 3$ 
\begin{equation*}
\int\limits_{0}^{2\pi }\int\limits_{\mathbb{R}^{n}}\frac{\left\vert
u(x,t)\right\vert ^{2}}{\left( x_{i}^{2}+x_{j}^{2}+x_{k}^{2}\right) ^{\delta
}}dxdt\leqslant C\left\Vert f\right\Vert _{2}^{2}.
\end{equation*}
This is equivalent to proving a uniform $L^{2}(\mathbb{R}^{n})$ boundedness
result for a family of singularized Hermite projection kernels.

As an application of the above estimate, we also prove the $\mathbb{R}^{9}$
collapsing variable type Strichartz estimate
\begin{equation*}
\int\limits_{0}^{2\pi }\int\limits_{\mathbb{R}^{3}}\left\vert u(\mathbf{x},
\mathbf{x},\mathbf{x},t)\right\vert ^{2}d\mathbf{x}dt\leqslant C\left\Vert
(-\triangle +\left\vert x\right\vert ^{2})f\right\Vert _{2}^{2}
\end{equation*}
where $\mathbf{x\in }\mathbb{R}^{3}$.
\end{abstract}

\maketitle

\section{Introduction}

\label{sec:intro}

In Bose-Einstein condensation (BEC), particles of integer spins
(\textquotedblleft Bosons\textquotedblright ) occupy a macroscopic quantum
state often called the \textquotedblleft condensation\textquotedblright . In
early lab experiments of BEC \cite{Anderson} \cite{Dyson}, the particles
were kept together by use of trapping potentials created by the effect of a
magnetic field on the particle spins. In principle, the magnetic field has a
complicated spatial structure. The interaction of the magnetic field with
the spin is conveniently modeled by a quadratic potential. This captures
salient features of the actual trap, especially the property that the
external potential rises at large distances. In later experiments, e.g., 
\cite{Stamper}, the trapping potential is produced by complicated laser
fields; but mathematically, one can still use a quadratic potential as a
simplified yet generic model. So the spin of the particle is removed in
modeling and the effect of a trap is included in the form of a quadratic
external potential. This physical background suggests that we study the Schr\"{o}dinger equation with quadratic potential

\begin{equation}
i\frac{\partial }{\partial t}u=-\triangle u+\left\vert x\right\vert ^{2}u
\text{ }in\text{ }\mathbb{R}^{n+1},  \label{quardtic equation}
\end{equation}
with initial data
\begin{equation*}
u(x,0)=f(x)\in L^{2}(\mathbb{R}^{n}).
\end{equation*}
Many aspects of equation \ref{quardtic equation} which came from the study
of the free Schr\"{o}dinger equation
\begin{equation}
i\frac{\partial }{\partial t}\phi =-\triangle \phi \text{ in }\mathbb{R}
^{n+1}  \label{FreeSchrodinger}
\end{equation}
have been studied by several authors. Its Strichartz estimates were proved
by Koch and Tataru \cite{BongioanniRogers}, Carles \cite{RemyCarles},
Nandakumarana and Ratnakuma \cite{Ratnakumar}. The well-posedness of its
nonlinear energy critical version with radial initial data was studied by
Killip, Visan and Zhang \cite{VisanAndZhang}. Bongioanni and Torrea, and
Bongioanni and Rogers, proved results on the pointwise convergence to the
initial data in \cite{Bongioanni} and \cite{BongioanniRogers}. Concerning
the Kato $\frac{1}{2}$-smoothing effect, Doi (and later Bongioanni and
Rogers) proved
\begin{equation}
\int\limits_{0}^{T}\int\limits_{\mathbb{R}^{n}}\frac{\left\vert
(I-(-\triangle +\left\vert x\right\vert ^{2}))^{\frac{1}{4}
}u(x,t)\right\vert ^{2}}{(1+\left\vert x\right\vert ^{2})^{\frac{1}{2}
+\varepsilon }}dxdt\leqslant C\left\Vert f\right\Vert _{2}^{2}
\label{Doi'sVegaEstimate}
\end{equation}
in \cite{Doi} (\cite{BongioanniRogers}), Robbiano and Zuily proved 
\begin{equation*}
\int\limits_{0}^{T}\int\limits_{\Omega }\left\vert \chi (x)(I-(-\triangle
+\left\vert x\right\vert ^{2}))^{\frac{1}{4}}u(x,t)\right\vert
^{2}dxdt\leqslant C\left\Vert f\right\Vert _{2}^{2}
\end{equation*}
for an external domain $\Omega $ and $\chi \in C_{0}^{\infty }(\overline{
\Omega })$ in \cite{Zuily}. However, both \cite{Doi} and \cite{Zuily} made
extensive use of pseudo-differential techniques which did not suffice to
prove the equation \ref{quardtic equation} counterpart to the Kato estimate 
\begin{equation}
\int\limits_{-\infty }^{\infty }\int\limits_{\mathbb{R}^{n}}\frac{
\left\vert \phi (x,t)\right\vert ^{2}}{\left\vert x\right\vert ^{2}}
dxdt\leqslant C\left\Vert \phi (\cdot ,0)\right\Vert _{2}^{2}
\label{Kato Estimate}
\end{equation}
in Kato and Yajima \cite{Kato}, or its generalization
\begin{equation}
\int\limits_{-\infty }^{\infty }\int\limits_{\mathbb{R}^{n}}\frac{
\left\vert \left\vert \nabla \right\vert ^{\alpha }\phi (x,t)\right\vert ^{2}
}{\left\vert x\right\vert ^{2-2\alpha }}dxdt\leqslant C\left\Vert \phi
(\cdot ,0)\right\Vert _{2}^{2},\text{ for }\alpha \in \lbrack 0,\frac{1}{2})
\label{KlainermanKato}
\end{equation}
in Kato and Yajima \cite{Kato}, and Ben-Artzi and Klainerman \cite
{Klainerman}, where $\phi $ is the solution to the free Schr\"{o}dinger
equation \ref{FreeSchrodinger}\ in the case $n\geqslant 3$.

\begin{remark}
Doi proved \ref{Doi'sVegaEstimate} type estimates in the case involving
variable coefficients. Bongioanni and Rogers proved an estimate similar to 
\ref{Doi'sVegaEstimate} for equation \ref{quardtic equation} using Hermite
functions. Their paper contains a series of results parallel to those in
Vega \cite{Vega}.
\end{remark}

\begin{remark}
The expository note \cite{ChenNote} gives extensions of estimate \ref
{KlainermanKato} to a class of dispersive equations, and simultaneously
arrives at the optimal constant for each $\alpha $ and $n$. The fact that $
\frac{\pi }{n-2}$ is the best constant achievable for the $\alpha =0$ case
is due to Simon \cite{Simon}.
\end{remark}

This paper aims to prove Kato type smoothing estimates similar to \ref{Kato
Estimate} when $n\geqslant 3$ for equation \ref{quardtic equation} without
using any pseudo-differential techniques.

In fact, we have

\begin{theorem}
\label{main} Let u be the solution to equation \ref{quardtic equation} in
the case $n\geqslant 3$, then for $\delta \in \lbrack 0,1]$ and $i\neq j\neq
k,$ one has the estimate
\begin{equation}
\int\limits_{0}^{2\pi }\int\limits_{\mathbb{R}^{n}}\frac{\left\vert
u(x,t)\right\vert ^{2}}{\left( x_{i}^{2}+x_{j}^{2}+x_{k}^{2}\right) ^{\delta
}}dxdt\leqslant C\left\Vert f\right\Vert _{2}^{2}.  \label{3d estimate}
\end{equation}
In particular, when $\delta =1,$ the above estimate implies the Kato
estimate \ref{Kato Estimate} for equation \ref{quardtic equation} in the
case $n\geqslant 3$ because of the trivial inequality $\left(
x_{i}^{2}+x_{j}^{2}+x_{k}^{2}\right) \leqslant \left\vert x\right\vert ^{2}.$
\end{theorem}

\begin{remark}
$u$ naturally has period $2\pi $ in the time variable $t$. We will show this
in section \ref{sec:HermiteFunctions}. This was also shown in Nandakumarana
and Ratnakuma \cite{Ratnakumar}.
\end{remark}

\begin{remark}
Without lose of generality, from here on out we assume $i=1,$ $j=2,$ $k=3$
for simplicity since the general case has an identical proof.
\end{remark}

\begin{remark}
We shall also point out that the Kato type estimate 
\begin{equation*}
\int\limits_{0}^{2\pi }\int\limits_{\mathbb{R}^{n}}\frac{\left\vert
\left\vert \nabla \right\vert ^{\alpha }u(x,t)\right\vert ^{2}}{\left\vert
x\right\vert ^{2-2\alpha }}dxdt\leqslant C\left\Vert f\right\Vert _{2}^{2},
\text{ for }\alpha \in \lbrack 0,\frac{1}{2})
\end{equation*}
which is the equation \ref{quardtic equation} counterpart to estimate \ref
{KlainermanKato} is still unproven.
\end{remark}

As an application of theorem \ref{main}, we have the following collapsing
variable type Strichartz estimate.

\begin{theorem}
\label{collapsing variables}Let $u(\mathbf{x}_{1},\mathbf{x}_{2},\mathbf{x}
_{3},t)$ solves equation \ref{quardtic equation} in the case $n=9$ with $
\mathbf{x}_{i}\in \mathbb{R}^{3}$, then one has the estimate 
\begin{equation}
\int\limits_{0}^{2\pi }\int\limits_{\mathbb{R}^{3}}\left\vert u(\mathbf{x},
\mathbf{x,x,}t)\right\vert ^{2}d\mathbf{x}dt\leqslant C\left\Vert
(-\triangle +\left\vert x\right\vert ^{2})f\right\Vert _{2}^{2}.
\label{collapsing variable estimate}
\end{equation}
\end{theorem}

\begin{remark}
The ordinary Strichartz estimate gives
\begin{equation*}
\left\Vert \phi \right\Vert _{L_{t}^{6}L_{x}^{6}}\leqslant C\left\Vert \phi
(\cdot ,0)\right\Vert _{\dot{H}^{\frac{2}{3}}}
\end{equation*}
for $\phi $ satisfying equation \ref{FreeSchrodinger} in the case $n=3.$
This leads us to consider estimate \ref{collapsing variable estimate}.
Estimates similar to \ref{collapsing variable estimate} were also considered
in Grillakis and Margetis \cite{Grillakis and Margetis}, and Klainerman and
Machedon \cite{KlainermanAndMachedon} in the setting of interacting Boson
systems.
\end{remark}

\begin{remark}
In Bongioanni and Torrea \cite{Bongioanni}, Bongioanni and Rogers \cite
{BongioanniRogers}, and Thangavelu \cite{Thangavelu1}, $\left\Vert
(-\triangle +\left\vert x\right\vert ^{2})^{\frac{s}{2}}f\right\Vert
_{2}^{2} $ is called the Hermite-Soblev $\mathcal{H}^{s}$ norm of $f.$ In
sections \ref{sec:Proof of the even estimate} and \ref{sec:collapsing
variables}, we will need the following lemma proved by Thangavelu concerning
Hermite-Soblev spaces in \cite{Thangavelu1}.
\end{remark}

\begin{lemma}
\label{HermiteSoblev}\cite{Thangavelu1} The operator $(I-\triangle )^{\frac{s
}{2}}(-\triangle +\left\vert x\right\vert ^{2})^{-\frac{s}{2}}$ is bounded
on $L^{2}(\mathbb{R}^{n})$ for $s\geqslant 0,$ or in other words
\begin{equation*}
\left\Vert (I-\triangle )^{\frac{s}{2}}f\right\Vert _{2}\leqslant
C\left\Vert (-\triangle +\left\vert x\right\vert ^{2})^{\frac{s}{2}
}f\right\Vert _{2},\text{ }s\geqslant 0.
\end{equation*}
\end{lemma}

Theorem \ref{main} will be deduced from the theorems below.

\begin{theorem}
\label{2d estimate} Let u be the solution to equation \ref{quardtic equation}
in the case $n\geqslant 2,$ then $\forall \delta \in \lbrack 0,1),$ one has
the estimate
\begin{equation*}
\int\limits_{0}^{2\pi }\int\limits_{\mathbb{R}^{n}}\frac{\left\vert
u(x,t)\right\vert ^{2}}{\left( x_{1}^{2}+x_{2}^{2}\right) ^{\delta }}
dxdt\leqslant C\left\Vert f\right\Vert _{2}^{2},
\end{equation*}
which implies estimate \ref{3d estimate} in the case $n\geqslant 3\ $and $
\delta \in \lbrack 0,1).$
\end{theorem}

\begin{theorem}
\label{1d odd estimate} Say $g(-x)=-$ $g(x)\in L^{2}(\mathbb{R})$, then we
have the equality
\begin{equation}
\int\limits_{0}^{2\pi }\int\limits_{\mathbb{R}}\frac{\left\vert
e^{-it(-\triangle +\left\vert x\right\vert ^{2})}g\right\vert ^{2}}{
\left\vert x\right\vert ^{2}}dxdt=4\pi \left\Vert g\right\Vert _{2}^{2}.
\label{odd estimate}
\end{equation}
In particular, estimate \ref{odd estimate} is equivalent to estimate \ref{3d
estimate} in the 3d radial case i.e. 
\begin{equation*}
\int\limits_{0}^{2\pi }\int\limits_{\mathbb{R}^{3}}\frac{\left\vert
e^{-it(-\triangle +\left\vert x\right\vert ^{2})}\psi \right\vert ^{2}}{
\left\vert x\right\vert ^{2}}dxdt=4\pi \left\Vert \psi \right\Vert _{L^{2}(
\mathbb{R}^{3})}^{2}
\end{equation*}
if $\psi $ is a $L^{2}(\mathbb{R}^{3})$ radial function.
\end{theorem}

\begin{remark}
There is an identity similar to equality \ref{odd estimate} for the free Schr\"{o}dinger equation \ref{FreeSchrodinger}. See the expository note \cite
{ChenNote}.
\end{remark}

\begin{theorem}
\label{3d even estimate}Say $d(\pm x_{1},\pm x_{2},\pm x_{3})=$ $
d(x_{1},x_{2},x_{3})\in L^{2}(\mathbb{R}^{3})$, then one has the estimate
\begin{equation*}
\int\limits_{0}^{2\pi }\int\limits_{\mathbb{R}^{3}}\frac{\left\vert
e^{-it(-\triangle +\left\vert x\right\vert ^{2})}d\right\vert ^{2}}{
\left\vert x\right\vert ^{2}}dxdt\leqslant C\left\Vert d\right\Vert _{2}^{2}.
\end{equation*}
\end{theorem}

Because we can write $f$ as a sum of its $x_{1}$-odd part and $x_{1}$-even
part by defining
\begin{equation*}
f_{odd}(x)=\frac{f(x_{1},x_{2},x_{3})-f(-x_{1},x_{2},x_{3})}{2}
\end{equation*}
and
\begin{equation*}
f_{even}(x)=\frac{f(-x_{1},x_{2},x_{3})+f(x_{1},x_{2},x_{3})}{2}.
\end{equation*}
So 
\begin{equation}
f(x)=f_{odd}(x)+f_{even,odd}(x)+f_{even,even,odd}(x)+f_{even,even,even}(x)
\label{odd-even decomposition}
\end{equation}
if we iterate the procedure three times. The linearity of equation \ref
{quardtic equation} and the fact that the terms in \ref{odd-even
decomposition} are all linear combinations of $f$ shows that estimate \ref
{3d estimate} in the case when $n=3$ and $\delta =1$ indeed follows from
theorems \ref{1d odd estimate} and \ref{3d even estimate}.

Moreover, theorems \ref{main} and \ref{2d estimate} are equivalent to the
following uniform $L^{2}(\mathbb{R}^{n})$ to $L^{2}(\mathbb{R}^{n})$
boundedness result for a family of singularized Hermite projection kernels.

\begin{theorem}
\label{Kernel estimate}For $n\geqslant 3$ and $\delta \in \lbrack 0,1],($ $
\delta \in \lbrack 0,1)$ when $n=2$ $)$, the singularized Hermite projection
kernels $\left\{ \frac{\Phi _{k}(x,y)}{\left\vert x\right\vert ^{\delta }}
\right\} _{k}$ map $L^{2}(\mathbb{R}^{n})$ to $L^{2}(\mathbb{R}^{n})$
uniformly where $\Phi _{k}$ is the usual Hermite projection kernel with
respect to the $k$-eigenspace defined in lemma \ref{ThangaveluLemma}, in
other words, there exists a $C>0$ depending only on $\delta $ and $n$ such
that
\begin{equation*}
\left\Vert \int_{\mathbb{R}^{n}}\frac{\Phi _{k}(\cdot ,y)}{\left\vert \cdot
\right\vert ^{\delta }}f(y)dy\right\Vert _{2}\leqslant C\left\Vert
f\right\Vert _{2}.
\end{equation*}
Moreover the more singular family $\left\{ \frac{\Phi _{k}(x,y)}{\left\vert
x\right\vert ^{\delta }\left\vert y\right\vert ^{\delta }}\right\} _{k}$
also maps $L^{2}(\mathbb{R}^{n})$ to $L^{2}(\mathbb{R}^{n})$ uniformly via
the standard $TT^{\ast \text{ }}$method.
\end{theorem}

Regularized Hermite projection kernels were studied in \cite{Bongioanni}, 
\cite{Ratnakumar}, \cite{Petrushev}, and \cite{Thangavelu}. But, to the best
of the author's knowledge, theorem \ref{Kernel estimate} might be the first
result on the singularized Hermite projection kernels.

\section{Some basics of Hermite functions and the proof of theorem \protect
\ref{2d estimate}}

\label{sec:HermiteFunctions}

To prove theorems \ref{main}, \ref{2d estimate}, \ref{1d odd estimate}, \ref
{3d even estimate} and \ref{Kernel estimate}, we will need the Hermite
functions and some of their properties. For more details outside of lemma 
\ref{Antiderivative} whose proof is provided in the appendix I, we refer the
reader to Thangavelu's monograph \cite{Thangavelu}.

\begin{definition}
\cite{Thangavelu}We define an $n$ dimensional Hermite function $\Phi
_{\alpha }(x)$ where $\alpha $ is an $n$-multiindex by 
\begin{equation*}
\Phi _{\alpha }(x)=\dprod\limits_{i=1}^{n}h_{\alpha _{i}}(x_{i}),
\end{equation*}
where $h_{k}$ are the one dimensional normalized Hermite functions defined by
\begin{equation*}
h_{k}(t)=\frac{(-1)^{k}}{(2^{k}k!\sqrt{\pi })^{\frac{1}{2}}}e^{\frac{t^{2}}{2
}}\frac{d^{k}}{dt^{k}}(e^{-t^{2}}),\text{ }t\in \mathbb{R}.
\end{equation*}
\end{definition}

Then we have the following well-known properties.

\begin{lemma}
\cite{Thangavelu}$\Phi _{\alpha }$ are the eigenfunctions of the Fourier
transform with eigenvalues $(-i)^{\left\vert \alpha \right\vert }$ i.e. 
\begin{equation*}
\widehat{\Phi _{\alpha }}(\xi )=(-i)^{\left\vert \alpha \right\vert }\Phi
_{\alpha }(\xi )
\end{equation*}
\end{lemma}

\begin{lemma}
\cite{Thangavelu} $\Phi _{\alpha }$ are also the eigenfunctions of the
Hermite operator $-\triangle +\left\vert x\right\vert ^{2}$ with eigenvalues 
$2\left\vert \alpha \right\vert +n.$ Moreover they form an orthonormal basis
of $L^{2}(\mathbb{R}^{n})$.
\end{lemma}

As this lemma states, we can write 
\begin{equation}
u(x,t)=\dsum\limits_{\alpha }e^{-i\lambda _{\alpha }t}a_{\alpha }\Phi
_{\alpha }(x),  \label{solutionformula}
\end{equation}
where $a_{\alpha }$ are the Fourier-Hermite coefficients
\begin{equation*}
a_{\alpha }=\int_{\mathbb{R}^{n}}f(x)\Phi _{\alpha }(x)dx,\text{ }
\end{equation*}
with convergence in $L^{2}(\mathbb{R}^{n})$, i.e. $u$ is naturally periodic $
2\pi $ in the time variable $t$ and we have
\begin{eqnarray*}
\int\limits_{0}^{2\pi }\left\vert u(x,t)\right\vert ^{2}dt
&=&\int\limits_{0}^{2\pi }\dsum\limits_{\alpha ,\beta }e^{-i(\lambda
_{\alpha }-\lambda _{\beta })t}a_{\alpha }\overline{a_{\beta }}\Phi _{\alpha
}(x)\Phi _{\beta }(x)dt \\
&=&2\pi \dsum\limits_{\substack{ \alpha ,\beta  \\ \lambda _{\alpha
}=\lambda _{\beta }}}a_{\alpha }\overline{a_{\beta }}\Phi _{\alpha }(x)\Phi
_{\beta }(x).
\end{eqnarray*}
But

\begin{equation*}
\int\limits_{\mathbb{R}}\dsum\limits_{\substack{ \alpha ,\beta  \\ \lambda
_{\alpha }=\lambda _{\beta }}}a_{\alpha }\overline{a_{\beta }}\Phi _{\alpha
}(x)\Phi _{\beta }(x)dx_{j}=\dsum\limits_{\substack{ \alpha ,\beta  \\ 
\lambda _{\alpha }=\lambda _{\beta }}}\delta _{a_{j}\beta _{j}}a_{\alpha }
\overline{a_{\beta }}\dprod\limits_{_{\substack{ i=1  \\ i\neq j}}
}^{n}h_{\alpha _{i}}(x_{i})\dprod\limits_{_{\substack{ i=1  \\ i\neq j}}
}^{n}h_{\beta _{i}}(x_{i}).
\end{equation*}
that is:
\begin{equation*}
\int\limits_{0}^{2\pi }\int\limits_{\mathbb{R}^{n+1}}\frac{\left\vert
u(x,t)\right\vert ^{2}}{\left\vert (x,x_{n+1})\right\vert ^{2\delta }}
dxdx_{n+1}dt\leqslant 2\pi \int\limits_{\mathbb{R}^{n}}\frac{1}{\left\vert
x\right\vert ^{2\delta }}\dsum\limits_{\substack{ \alpha ,\beta  \\ \lambda
_{\alpha }=\lambda _{\beta }}}\delta _{a_{j}\beta _{j}}a_{\alpha }\overline{
a_{\beta }}\dprod\limits_{_{i=1}}^{n}h_{\alpha
_{i}}(x_{i})\dprod\limits_{_{i=1}}^{n}h_{\beta _{i}}(x_{i})dx
\end{equation*}
hence the estimate
\begin{equation}
\int\limits_{0}^{2\pi }\int\limits_{\mathbb{R}^{2}}\frac{\left\vert
u(x,t)\right\vert ^{2}}{\left\vert x\right\vert ^{2\delta }}dxdt\leqslant
C\left\Vert f\right\Vert _{2}^{2}  \label{2d estimate 1}
\end{equation}
for $\delta \in \lbrack 0,1)$ implies theorem \ref{2d estimate} and we only
need to prove theorem \ref{main} in the case $n=3$ and $\delta =1.$ We can
now prove estimate $\ref{2d estimate 1}$ by the following lemma.

\begin{lemma}
\label{ThangaveluLemma}\cite{Thangavelu} Let $P_{k}$ be the Hermite
projector corresponding to the $k$-eigenspace with kernel
\begin{equation*}
\Phi _{k}(x,y)=\dsum\limits_{\left\vert \alpha \right\vert =k}\Phi _{\alpha
}(x)\Phi _{\alpha }(y)
\end{equation*}
then there is a constant $C\geqslant 0$ independent of $k$ and $x$ such that
\begin{equation}
\left\vert \Phi _{k}(x,x)\right\vert \leqslant Ck^{\frac{n}{2}-1}.
\label{Thangavelu estimate}
\end{equation}
\end{lemma}

Therefore
\begin{eqnarray}
&&\int\limits_{0}^{2\pi }\int\limits_{\mathbb{R}^{2}}\frac{\left\vert
u(x,t)\right\vert ^{2}}{\left\vert x\right\vert ^{2\delta }}dxdt
\label{pde to kernel} \\
&=&2\pi \dsum\limits_{k}\int\limits_{\mathbb{R}^{2}}\frac{\left\vert
P_{k}f\right\vert ^{2}}{\left\vert x\right\vert ^{2\delta }}dx  \notag \\
&\leqslant &2\pi \dsum\limits_{k}\int\limits_{\mathbb{R}^{2}-\mathbb{D}
}\left\vert P_{k}f\right\vert ^{2}dx+2\pi \dsum\limits_{k}\int_{\mathbb{D}}
\frac{\left\vert P_{k}f\right\vert ^{2}}{\left\vert x\right\vert ^{2\delta }}
dx  \notag \\
&\leqslant &2\pi \left\Vert f\right\Vert _{2}^{2}+2\pi \dsum\limits_{k}\int_{
\mathbb{D}}\frac{1}{\left\vert x\right\vert ^{2\delta }}\bigg(\dsum\limits 
_{\substack{ \alpha ,\beta  \\ \left\vert \alpha \right\vert =\left\vert
\beta \right\vert =k}}\left\vert a_{\alpha }\right\vert ^{2}\left\vert 
\overline{a_{\beta }}\right\vert ^{2}\bigg)^{\frac{1}{2}}\bigg(
\dsum\limits_{ _{\substack{ \alpha ,\beta  \\ \left\vert \alpha \right\vert
=\left\vert \beta \right\vert =k}}}\left\vert \Phi _{\alpha }(x)\Phi _{\beta
}(x)\right\vert ^{2}\bigg)^{\frac{1}{2}}dx  \notag \\
&\leqslant &2\pi \left\Vert f\right\Vert _{2}^{2}+C\dsum\limits_{k}\bigg(
\left\Vert P_{k}f\right\Vert _{2}^{2}\int\limits_{0}^{2\pi
}\int\limits_{0}^{1}\frac{1}{r^{2\delta -1}}drd\theta \bigg)  \notag \\
&\leqslant &C\left\Vert f\right\Vert _{2}^{2}.  \notag
\end{eqnarray}

\begin{remark}
In the above computation, we have also proved the Morawetz inequality
\begin{equation*}
\sup_{x\in \mathbb{R}^{2}}\int\limits_{0}^{2\pi }\left\vert
u(x,t)\right\vert ^{2}dt\leqslant C\left\Vert f\right\Vert _{2}^{2}
\end{equation*}
which is identical to the well-known version for the free Schr\"{o}dinger
equation \ref{FreeSchrodinger} in \cite{TaoMora}. This is another Kato type
smoothing estimate.

\begin{remark}
Estimate \ref{Thangavelu estimate} is also the key ingredient to prove the
regularized Hermite projection kernel estimates in \cite{Petrushev}. But it
does not yield theorem \ref{main} in the case when $n=3$ and $\delta =1.$
Lemma \ref{Antiderivative} will introduce a new tool for that purpose.
\end{remark}
\end{remark}

As
\begin{equation*}
\int\limits_{0}^{2\pi }\int\limits_{\mathbb{R}^{3}}\frac{\left\vert
u(x,t)\right\vert ^{2}}{\left\vert x\right\vert ^{2\delta }}dxdt=2\pi
\dsum\limits_{k}\int\limits_{\mathbb{R}^{3}}\frac{\left\vert
P_{k}f\right\vert ^{2}}{\left\vert x\right\vert ^{2\delta }}dx,
\end{equation*}
theorem \ref{Kernel estimate} implies theorems \ref{main} and \ref{2d
estimate}. However, the fact that $e^{-i(2k+n)t}P_{k}f$ satisfies equation 
\ref{quardtic equation} with $u(x,0)=P_{k}f(x)$ shows that theorems \ref
{main} and \ref{2d estimate} also imply theorem \ref{Kernel estimate}.

We are left with the proofs of theorems \ref{1d odd estimate} and \ref{3d
even estimate} which will need the following tool.

\begin{lemma}
\label{Antiderivative}We define the "antiderivatives" of the 1-d Hermite
functions to be
\begin{equation*}
X_{2k+1}(x)=\int_{-\infty }^{x}h_{2k+1}(t)dt
\end{equation*}
and
\begin{equation*}
X_{2k}(x)=\int_{-\infty }^{x}sign(t)h_{2k}(t)dt
\end{equation*}
which are by definition absolutely continuous. Moreover
\begin{equation}
\int_{\mathbb{R}}(X_{2k+1}(x))^{2}dx=2  \label{odd computation}
\end{equation}
and
\begin{eqnarray}
\int_{\mathbb{R}}(X_{2k}(x))^{2}dx &=&2(-1+\sqrt{2}\dsum\limits_{i=0}^{k}
\binom{\frac{1}{2}}{i})  \label{even computation} \\
&\leqslant &3  \notag
\end{eqnarray}
where
\begin{equation*}
\dsum\limits_{i=0}^{\infty }\binom{\frac{1}{2}}{i}=\sqrt{2}
\end{equation*}
i.e. $\lim_{k\rightarrow \infty }\left\Vert X_{2k}\right\Vert _{2}^{2}$ $=2$
and $X_{k}\in H^{1}(\mathbb{R}).$
\end{lemma}

To the best of our knowledge, lemma \ref{Antiderivative} is new. The proof,
which is a direct computation, is provided in the appendix I\ for
completion. Now we can give the proofs of theorems \ref{1d odd estimate} and 
\ref{3d even estimate}.

\section{Proof of theorem \protect\ref{1d odd estimate}}

\label{sec:Proof of the odd estimate}We only need to prove 
\begin{equation*}
\int\limits_{\mathbb{R}}\frac{\left\vert P_{2k+1}g\right\vert ^{2}}{
\left\vert x\right\vert ^{2}}dx=2\left\vert a_{2k+1}\right\vert ^{2}
\end{equation*}

because $g(-x)=-$ $g(x)$ implies $a_{2k}=0$, $\forall k$. In fact since $
h_{2k}(x)$ is even, we have
\begin{eqnarray*}
a_{2k} &=&\int_{\mathbb{R}}g(x)h_{2k}(x)dx \\
&=&0.
\end{eqnarray*}

One notices that 
\begin{equation*}
h_{2k+1}(\xi )=\frac{d}{dx}X_{2k+1}(\xi )
\end{equation*}
i.e.
\begin{equation*}
(-i)^{-(2k+1)}h_{2k+1}(x)=x\check{X}_{2k+1}(x)
\end{equation*}
hence 
\begin{eqnarray*}
\int\limits_{\mathbb{R}}\frac{\left\vert P_{2k+1}g\right\vert ^{2}}{
\left\vert x\right\vert ^{2}}dx &=&\left\vert a_{2k+1}\right\vert
^{2}\int\limits_{\mathbb{R}}\left\vert \check{X}_{2k+1}(x)\right\vert ^{2}dx
\\
&=&\left\vert a_{2k+1}\right\vert ^{2}\int\limits_{\mathbb{R}}\left\vert
X_{2k+1}(x)\right\vert ^{2}dx \\
&=&2\left\vert a_{2k+1}\right\vert ^{2}
\end{eqnarray*}
via equality \ref{odd computation}. Whence we have deduced theorem \ref{1d
odd estimate}.

\section{Proof of theorem \protect\ref{3d even estimate}}

\label{sec:Proof of the even estimate}It suffices to prove that there exists
a $C>0$ independent of $k$ s.t.
\begin{equation*}
\int\limits_{\mathbb{R}^{3}}\frac{\left\vert P_{k}d\right\vert ^{2}}{
\left\vert x\right\vert ^{2}}dx\leqslant C\left\Vert P_{k}d\right\Vert
_{2}^{2}.
\end{equation*}

Throughout this section, we will assume $k\neq 0$. In the case when $k=0$, $
P_{k}d$ has only one term
\begin{eqnarray*}
P_{0}d &=&a_{0}h_{0}(x_{1})h_{0}(x_{2})h_{0}(x_{3}) \\
&=&a_{0}\left( \sqrt{\pi }\right) ^{-\frac{3}{2}}e^{-\frac{\left\vert
x\right\vert ^{2}}{2}}
\end{eqnarray*}
and hence is a 3d radial function, and we dealt with this situation in
theorem \ref{1d odd estimate}. In fact, it is easy to compute that
\begin{equation*}
\int\limits_{\mathbb{R}^{3}}\frac{\left\vert P_{0}d\right\vert ^{2}}{
\left\vert x\right\vert ^{2}}dx=2\left\vert a_{0}\right\vert ^{2}
\end{equation*}
which matches theorem \ref{1d odd estimate}.

We write 
\begin{eqnarray*}
P_{k}d(x) &=&\dsum\limits_{_{\substack{ \alpha  \\ \alpha _{1},\text{ }
\alpha _{2},\text{ }\alpha _{3}\text{ are even}  \\ \left\vert \alpha
\right\vert =k}}}a_{\alpha }h_{\alpha _{1}}(x_{1})h_{\alpha
_{2}}(x_{2})h_{\alpha _{3}}(x_{3}) \\
&=&\sum_{\alpha \in I}+\sum_{\alpha \in II}+\sum_{\alpha \in III}
\end{eqnarray*}
where 
\begin{equation*}
I=\{\alpha :\left\vert \alpha \right\vert =k,\text{ }\alpha _{1},\text{ }
\alpha _{2},\text{ }\alpha _{3}\text{ are even, and }\alpha _{1}\geqslant 
\frac{\alpha _{2}+\text{ }\alpha _{3}}{2}\}
\end{equation*}
\begin{equation*}
II=\{\alpha :\left\vert \alpha \right\vert =k,\text{ }\alpha _{1},\text{ }
\alpha _{2},\text{ }\alpha _{3}\text{ are even, and }\alpha _{2}\geqslant 
\frac{\alpha _{1}+\text{ }\alpha _{3}}{2}\}
\end{equation*}
\begin{equation*}
III=\{\alpha :\left\vert \alpha \right\vert =k,\text{ }\alpha _{1},\text{ }
\alpha _{2},\text{ }\alpha _{3}\text{ are even, and }\alpha _{3}\geqslant 
\frac{\alpha _{1}+\text{ }\alpha _{2}}{2}\}.
\end{equation*}

\begin{remark}
Suppose we have $\alpha $ s.t. $\alpha _{1}<\frac{\alpha _{2}+\text{ }\alpha
_{3}}{2},\alpha _{2}<\frac{\alpha _{1}+\text{ }\alpha _{3}}{2},$and $\alpha
_{3}<\frac{\alpha _{2}+\text{ }\alpha _{1}}{2}$, then $\alpha _{1}+\alpha
_{2}+\alpha _{3}<\alpha _{1}+\alpha _{2}+\alpha _{3}$ which is a
contradiction$.$ So $I,$ $II,$ and $III$ covers all cases. In some cases, $
I, $ $II,III$ do not intersect trivially. In these cases, we just count the
crossing terms once in one proper set. Moreover $a_{\alpha }=0,$ $\forall
\alpha $ with one odd index due to $d(\pm x_{1},\pm x_{2},\pm x_{3})=$ $
d(x_{1},x_{2},x_{3})$, in fact since $h_{\alpha _{1}}(x_{1})$ is odd if $
\alpha _{1}$ is odd, we have
\begin{eqnarray*}
a_{\alpha } &=&\int_{\mathbb{R}^{2}}dx_{2}dx_{3}\int_{\mathbb{R}
}dx_{1}d(x)h_{\alpha _{1}}(x_{1})h_{\alpha _{2}}(x_{2})h_{\alpha _{3}}(x_{3})
\\
&=&\int_{\mathbb{R}^{2}}dx_{2}dx_{3}h_{\alpha _{2}}(x_{2})h_{\alpha
_{3}}(x_{3})\cdot 0 \\
&=&0
\end{eqnarray*}
\end{remark}

\bigskip So it is enough to prove
\begin{equation}
\int\limits_{\mathbb{R}^{3}}\frac{\left\vert \sum_{\alpha \in I}\right\vert
^{2}}{\left\vert x\right\vert ^{2}}dx\leqslant C\sum_{\alpha \in
I}\left\vert a_{\alpha }\right\vert ^{2},  \label{Estimate I}
\end{equation}
\begin{equation}
\int\limits_{\mathbb{R}^{3}}\frac{\left\vert \sum_{\alpha \in
II}\right\vert ^{2}}{\left\vert x\right\vert ^{2}}dx\leqslant C\sum_{\alpha
\in II}\left\vert a_{\alpha }\right\vert ^{2},  \label{Estimate II}
\end{equation}
and
\begin{equation}
\int\limits_{\mathbb{R}^{3}}\frac{\left\vert \sum_{\alpha \in
III}\right\vert ^{2}}{\left\vert x\right\vert ^{2}}dx\leqslant C\sum_{\alpha
\in III}\left\vert a_{\alpha }\right\vert ^{2}.  \label{Estimate III}
\end{equation}

In the following, we will only prove estimate \ref{Estimate I}, and the
proofs of estimates \ref{Estimate II} and \ref{Estimate III} will be
similar. To be more specific, we will use $\alpha _{1}$ and $x_{1}$ for
estimate \ref{Estimate I}, $\alpha _{2}$ and $x_{2}$ for estimate \ref
{Estimate II}, $\alpha _{3}$ and $x_{3}$ for estimate \ref{Estimate III}.

Define
\begin{equation*}
u_{k,I}(\xi )=\sum_{\alpha \in I}a_{\alpha }X_{\alpha _{1}}(\xi
_{1})h_{\alpha _{2}}(\xi _{2})h_{\alpha _{3}}(\xi _{3})
\end{equation*}
then
\begin{eqnarray*}
&&\left\Vert u_{k,I}\right\Vert _{2}^{2} \\
&=&\dsum\limits_{\alpha ,\beta \in I}\int_{\mathbb{R}}d\xi _{1}\int_{\mathbb{
R}^{2}}d\xi _{2}d\xi _{3}a_{\alpha }\overline{a_{\beta }}X_{\alpha _{1}}(\xi
_{1})h_{\alpha _{2}}(\xi _{2})h_{\alpha _{3}}(\xi _{3})X_{\beta _{1}}(\xi
_{1})h_{\beta _{2}}(\xi _{2})h_{\beta _{3}}(\xi _{3}) \\
&=&\dsum\limits_{\alpha \in I}\left\vert a_{\alpha }\right\vert ^{2}\int_{
\mathbb{R}}(X_{\alpha _{1}}(\xi _{1}))^{2}d\xi _{1} \\
&\leqslant &3\sum_{\alpha \in I}\left\vert a_{\alpha }\right\vert ^{2},
\end{eqnarray*}
via formula \ref{even computation}.

Moreover, 
\begin{equation*}
\sum_{\alpha \in I}a_{\alpha }h_{\alpha _{1}}(\xi _{1})h_{\alpha _{2}}(\xi
_{2})h_{\alpha _{3}}(\xi _{3})=sign(\xi _{1})\frac{\partial }{\partial \xi
_{1}}u_{k,I}(\xi ),\text{ }\xi _{1}\neq 0
\end{equation*}
yields 
\begin{eqnarray*}
&&(-i)^{-k}\sum_{\alpha \in I}a_{\alpha }h_{\alpha 1}(x_{1})h_{\alpha
_{2}}(x_{2})h_{\alpha _{3}}(x_{3}) \\
&=&H(t_{1}\check{u}_{k,I}(t))(x) \\
&=&x_{1}H(\check{u}_{k,I})(x)-\int\limits_{-\infty }^{\infty }\check{u}
_{k,even}(t,x_{2},x_{3})dt
\end{eqnarray*}
where $H$ is the Hilbert transform only with respect to the the first
variable.

Hence we have
\begin{eqnarray*}
\int\limits_{\mathbb{R}^{3}}\frac{\left\vert \sum_{\alpha \in I}\right\vert
^{2}}{\left\vert x\right\vert ^{2}}dx &=&\int\limits_{\mathbb{R}^{3}}\frac{
\left\vert x_{1}H(\check{u}_{k,I})(x)-\int\limits_{-\infty }^{\infty }
\check{u}_{k,I}(t,x_{2},x_{3})dt\right\vert ^{2}}{\left\vert x\right\vert
^{2}}dx \\
&\leqslant &2\int\limits_{\mathbb{R}^{3}}\frac{\left\vert x_{1}H(\check{u}
_{k,I})(x)\right\vert ^{2}}{\left\vert x\right\vert ^{2}}dx+2\int\limits_{
\mathbb{R}^{3}}\frac{\left\vert \int\limits_{-\infty }^{\infty }\check{u}
_{k,I}(t,x_{2},x_{3})dt\right\vert ^{2}}{\left\vert x\right\vert ^{2}}
dx_{1}dx_{2}dx_{3} \\
&\leqslant &6\sum_{\alpha \in I}\left\vert a_{\alpha }\right\vert ^{2}+2\pi
\int\limits_{\mathbb{R}^{2}}\left\vert \int\limits_{-\infty }^{\infty }
\check{u}_{k,I}(t,x_{2},x_{3})dt\right\vert ^{2}\frac{dx_{2}dx_{3}}{\sqrt{
x_{2}^{2}+x_{3}^{2}}} \\
&=&6\sum_{\alpha \in I}\left\vert a_{\alpha }\right\vert ^{2}+2\pi
(MainTerm_{I})
\end{eqnarray*}
where
\begin{eqnarray*}
MainTerm_{I} &=&\int\limits_{\mathbb{R}^{2}}\left\vert
\int\limits_{-\infty }^{\infty }\check{u}_{k,I}(t,x_{2},x_{3})dt\right\vert
^{2}\frac{dx_{2}dx_{3}}{\sqrt{x_{2}^{2}+x_{3}^{2}}} \\
&=&\int\limits_{\mathbb{R}^{2}}\left\vert \sum_{\alpha \in
I}\int\limits_{-\infty }^{\infty }dt\int\limits_{\mathbb{R}^{3}}e^{it\xi
_{1}}e^{ix_{2}\xi _{2}}e^{ix_{3}\xi _{3}}a_{\alpha }X_{\alpha _{1}}(\xi
_{1})h_{\alpha _{2}}(\xi _{2})h_{\alpha _{3}}(\xi _{3})d\xi \right\vert ^{2}
\frac{dx_{2}dx_{3}}{\sqrt{x_{2}^{2}+x_{3}^{2}}} \\
&=&\int\limits_{\mathbb{R}^{2}}\left\vert \sum_{\alpha \in I}a_{\alpha
}X_{\alpha _{1}}(0)(-i)^{\alpha _{2}+\alpha _{3}}h_{\alpha
_{2}}(x_{2})h_{\alpha _{3}}(x_{3})\right\vert ^{2}\frac{dx_{2}dx_{3}}{\sqrt{
x_{2}^{2}+x_{3}^{2}}} \\
&\leqslant &C\int\limits_{\mathbb{R}^{2}}\left\vert \left\vert \nabla
\right\vert ^{\frac{1}{2}}(\sum_{\alpha \in I}a_{\alpha }X_{\alpha
_{1}}(0)(-i)^{\alpha _{2}+\alpha _{3}}h_{\alpha _{2}}(x_{2})h_{\alpha
_{3}}(x_{3}))\right\vert ^{2}dx_{2}dx_{3}\text{ }\left( \text{Hardy's
inequality}\right) \\
&\leqslant &C\int\limits_{\mathbb{R}^{2}}\left\vert (-\triangle +\left\vert
(x_{2},x_{3})\right\vert ^{2})^{\frac{1}{4}}{\LARGE (}\sum_{\alpha \in
I}a_{\alpha }X_{\alpha _{1}}(0)(-i)^{\alpha _{2}+\alpha _{3}}h_{\alpha
_{2}}(x_{2})h_{\alpha _{3}}(x_{3}){\LARGE )}\right\vert ^{2}dx_{2}dx_{3}
\text{ }\left( \text{Lemma \ref{HermiteSoblev}}\right) \\
&=&C\sum_{\alpha \in I}\left\vert a_{\alpha }\right\vert ^{2}(X_{\alpha
_{1}}(0))^{2}(2\alpha _{2}+2\alpha _{3}+2)^{\frac{1}{2}}.
\end{eqnarray*}

However, from Feldheim \cite{Feldheim} and Busbridge \cite{BusBridge}, we
know that given $\alpha _{1}$ even

\begin{eqnarray*}
(X_{\alpha _{1}}(0))^{2} &=&\frac{1}{4}(\int\limits_{-\infty }^{\infty
}h_{\alpha _{1}}(t)dt)^{2} \\
&=&\frac{\sqrt{2}}{4}\frac{2^{2\alpha _{1}}(\Gamma (\frac{1}{2}\alpha _{1}+
\frac{1}{2}))^{2}}{2^{\alpha _{1}}\alpha _{1}!\sqrt{\pi }}\text{ } \\
&\leqslant &C\frac{1}{(\alpha _{1})^{\frac{1}{2}}}
\end{eqnarray*}
by Stirling's formula. The above inequality shows
\begin{equation*}
(X_{\alpha _{1}}(0))^{2}(\alpha _{2}+\alpha _{3}+1)^{\frac{1}{2}}\leqslant C
\frac{(2\alpha _{1}+1)^{\frac{1}{2}}}{(\alpha _{1})^{\frac{1}{2}}}\leqslant C
\end{equation*}
for $\alpha _{1}\geqslant \frac{\alpha _{2}+\text{ }\alpha _{3}}{2}$ and $
\alpha _{1}\neq 0,$ or in other words, for $\alpha \in I$ and $k\neq 0$.

So
\begin{equation*}
Mainterm_{I}\leqslant C\sum_{\alpha \in I}\left\vert a_{\alpha }\right\vert
^{2}
\end{equation*}
i.e.
\begin{equation*}
\int\limits_{\mathbb{R}^{3}}\frac{\left\vert \sum_{\alpha \in I}\right\vert
^{2}}{\left\vert x\right\vert ^{2}}dx\leqslant C\sum_{\alpha \in
I}\left\vert a_{\alpha }\right\vert ^{2}
\end{equation*}

\begin{remark}
If we apply this procedure to the case when $n=2$ and $\delta
=1,Mainterm_{I} $ will have $\left\vert x_{2}\right\vert ^{-1}$ as a
singularity which forces $MainTerm_{I}$ to be $\infty $ whenever there is
some $a_{a}\neq 0$. But this procedure does also prove estimate \ref{2d
estimate 1} when $\delta <1$ and hence theorem \ref{2d estimate}.
\end{remark}

\section{An Application of Theorem \protect\ref{main} / Proof of Theorem 
\protect\ref{collapsing variables}}

\label{sec:collapsing variables}To obtain theorem \ref{collapsing variables}, aside from theorem \ref{main} and lemma \ref{HermiteSoblev}, an
interaction Morawetz inequality is needed.

\subsection{Morawetz inequality}

As in \cite{TaoMora}, define 
\begin{equation*}
T_{00}=\left\vert u\right\vert ^{2}
\end{equation*}
\begin{equation*}
T_{0j}=T_{j0}=2\func{Im}\frac{\partial u}{\partial x_{j}}\overline{u}
\end{equation*}
and
\begin{equation*}
T_{jk}=T_{kj}=4\func{Re}(u_{k}\overline{u_{j}})-\delta _{jk}\triangle
(\left\vert u\right\vert ^{2})
\end{equation*}
where $j$, $k$ mean summation from $1$ to $n$. Then a direct computation
shows that
\begin{eqnarray*}
\partial _{t}T_{00}+\partial _{j}T_{0j} &=&0, \\
\partial _{t}T_{k0}+\partial _{j}T_{kj} &=&-2V_{k}\left\vert u\right\vert
^{2},
\end{eqnarray*}
for the equation 
\begin{equation*}
iu_{t}=-\triangle u+Vu.
\end{equation*}

Hence we have 
\begin{eqnarray*}
\partial _{t}M_{0}^{a}(t) &=&4\int_{\mathbb{R}^{n}}a_{kj}\func{Re}(u_{k}
\overline{u_{j}})dx-\int_{\mathbb{R}^{n}}\triangle a\triangle (\left\vert
u\right\vert ^{2})dx-2\int_{\mathbb{R}^{n}}a_{k}V_{k}\left\vert u\right\vert
^{2}dx. \\
&=&4\int_{\mathbb{R}^{n}}a_{kj}\func{Re}(u_{k}\overline{u_{j}})dx-\int_{
\mathbb{R}^{n}}\triangle a\triangle (\left\vert u\right\vert ^{2})dx+2\int_{
\mathbb{R}^{n}}a(x)V_{kk}\left\vert u\right\vert ^{2}dx \\
&&+2\int_{\mathbb{R}^{n}}a(x)(\left\vert u\right\vert ^{2})_{k}V_{k}dx
\end{eqnarray*}
if we define
\begin{equation*}
M_{0}^{a}(t)=\int_{\mathbb{R}^{n}}a_{k}(x)T_{0k}(t,x)dx
\end{equation*}
to be the Morawetz action corresponding to a suitable $a(x)$ which will be
chosen momentarily.

Therefore, for equation \ref{quardtic equation} in the case $n=9$, we in
fact have
\begin{eqnarray*}
&&\int\limits_{0}^{2\pi }\int_{\mathbb{R}^{9}}(\triangle ^{2}a)\left\vert
u\right\vert ^{2}dxdt \\
&=&4\int\limits_{0}^{2\pi }\int_{\mathbb{R}^{9}}a_{kj}(x)\func{Re}(u_{k}
\overline{u_{j}})dxdt \\
&&+36\int\limits_{0}^{2\pi }\int_{\mathbb{R}^{9}}a(x)\left\vert
u\right\vert ^{2}dxdt \\
&&+4\int\limits_{0}^{2\pi }\int_{\mathbb{R}^{9}}a(x)V_{k}\func{Re}(u
\overline{u_{k}})dx
\end{eqnarray*}
due to the facts that 
\begin{equation*}
\int\limits_{0}^{2\pi }\partial _{t}M_{0}^{a}(t)dt=0
\end{equation*}
and $V_{kk}=18$.

The $a(x)$ we are going to pick is not non-strictly convex as in the usual
cases in \cite{TaoMora}, but the follwing computation will help to simplify
the technical problems arised from that:
\begin{eqnarray*}
&&2\int_{\mathbb{R}^{9}}a_{kj}\func{Re}(u_{k}\overline{u_{j}})dx \\
&=&\int_{\mathbb{R}^{9}}a_{kj}(u_{k}\overline{u_{j}}+u_{j}\overline{u_{k}})dx
\\
&=&-\int_{\mathbb{R}^{9}}a_{k}(u_{kj}\overline{u_{j}}+u_{k}\overline{u_{jj}}
+u_{jj}\overline{u_{k}}+u_{j}\overline{u_{kj}})dx \\
&=&-\int_{\mathbb{R}^{9}}a_{k}(u_{kj}\overline{u_{j}}+u_{j}\overline{u_{kj}}
)dx-\int_{\mathbb{R}^{9}}a_{k}(u_{k}\overline{\triangle u}+\overline{u_{k}}
\triangle u)dx \\
&=&-\int_{\mathbb{R}^{9}}a_{k}(\left\vert u_{j}\right\vert ^{2})_{k}dx+\int_{
\mathbb{R}^{9}}a(u_{kk}\overline{\triangle u}+u_{k}\overline{\triangle u_{k}}
+\overline{u_{kk}}\triangle u+\overline{u_{k}}\triangle u_{k})dx \\
&=&\int_{\mathbb{R}^{9}}a_{kk}\left\vert u_{j}\right\vert ^{2}dx+2\int_{
\mathbb{R}^{9}}a\left\vert \triangle u\right\vert ^{2}dx+\int_{\mathbb{R}
^{9}}a(u_{k}\overline{\triangle u_{k}}+\overline{u_{k}}\triangle u_{k})dx \\
&=&\int_{\mathbb{R}^{9}}\triangle a\left\vert \nabla u\right\vert
^{2}dx+2\int_{\mathbb{R}^{9}}a\left\vert \triangle u\right\vert ^{2}dx+\int_{
\mathbb{R}^{9}}a(\triangle \left\vert u_{k}\right\vert ^{2}-2\left\vert
\nabla u_{k}\right\vert ^{2})dx \\
&=&2\int_{\mathbb{R}^{9}}\triangle a\left\vert \nabla u\right\vert
^{2}dx+2\int_{\mathbb{R}^{9}}a\left\vert \triangle u\right\vert
^{2}dx-2\int_{\mathbb{R}^{9}}a\left\vert \nabla ^{2}u\right\vert ^{2}dx.
\end{eqnarray*}

So
\begin{eqnarray}
\int\limits_{0}^{2\pi }\int_{\mathbb{R}^{9}}(\triangle ^{2}a)\left\vert
u\right\vert ^{2}dxdt &=&4\int\limits_{0}^{2\pi }\int_{\mathbb{R}
^{9}}\triangle a\left\vert \nabla u\right\vert ^{2}dxdt  \label{Morawetz} \\
&&+4\int\limits_{0}^{2\pi }\int_{\mathbb{R}^{9}}a(x)\left\vert \triangle
u\right\vert ^{2}dxdt  \notag \\
&&-4\int\limits_{0}^{2\pi }\int_{\mathbb{R}^{9}}a(x)\left\vert \nabla
^{2}u\right\vert ^{2}dxdt  \notag \\
&&+36\int\limits_{0}^{2\pi }\int_{\mathbb{R}^{9}}a(x)\left\vert
u\right\vert ^{2}dxdt  \notag \\
&&+4\int\limits_{0}^{2\pi }\int_{\mathbb{R}^{9}}a(x)V_{k}\func{Re}(u
\overline{u_{k}})dxdt  \notag
\end{eqnarray}
If we select 
\begin{equation}
a(\mathbf{x}_{1},\mathbf{x}_{2},\mathbf{x}_{3})=C\frac{1}{\left\vert \mathbf{x}_{1}-\mathbf{x}_{2}\right\vert ^{2}+\left\vert \mathbf{x}_{1}-\mathbf{x}
_{3}\right\vert ^{2}+\left\vert \mathbf{x}_{2}-\mathbf{x}_{3}\right\vert ^{2}
}  \label{a}
\end{equation}
where $C$ is a suitable positive constant, then
\begin{equation*}
\triangle ^{2}a=\delta (\mathbf{x}_{1}-\mathbf{x}_{2})\delta (\mathbf{x}_{2}-
\mathbf{x}_{3}),
\end{equation*}
\begin{equation*}
\triangle a(\mathbf{x}_{1},\mathbf{x}_{2},\mathbf{x}_{3})=-C\frac{1}{
(\left\vert \mathbf{x}_{1}-\mathbf{x}_{2}\right\vert ^{2}+\left\vert \mathbf{x}_{1}-\mathbf{x}_{3}\right\vert ^{2}+\left\vert \mathbf{x}_{2}-\mathbf{x}
_{3}\right\vert ^{2})^{2}}<0,
\end{equation*}
and relation \ref{Morawetz} reads 
\begin{eqnarray*}
\int\limits_{0}^{2\pi }\int_{\mathbb{R}^{3}}\left\vert u(\mathbf{x},\mathbf{x},\mathbf{x,}t)\right\vert ^{2}d\mathbf{x}dt &\leqslant
&4\int\limits_{0}^{2\pi }\int_{\mathbb{R}^{9}}a\left\vert \triangle
u\right\vert ^{2}dxdt+36\int\limits_{0}^{2\pi }\int_{\mathbb{R}
^{9}}a(x)\left\vert u\right\vert ^{2}dxdt \\
&&+4\int\limits_{0}^{2\pi }\int_{\mathbb{R}^{9}}a(x)V_{k}\func{Re}(u
\overline{u_{k}})dxdt \\
&=&4A+36B+4D
\end{eqnarray*}

\begin{remark}
Formula \ref{a} is from Klainerman and Machedon's private communication\cite
{K&M private}. Thanks to Machedon for sharing this computation.
\end{remark}

To prove estimate \ref{collapsing variable estimate}, it will suffice to
show that $A,B,$and $D$ are majorized by $\left\Vert (-\triangle +\left\vert
x\right\vert ^{2})f\right\Vert _{2}^{2}.$

\subsection{Estimates for $A,B$, and $D$}

\begin{eqnarray*}
A &=&\int\limits_{0}^{2\pi }\int_{\mathbb{R}^{9}}a\left\vert \triangle
u\right\vert ^{2}dxdt \\
&=&C\int\limits_{0}^{2\pi }\int_{\mathbb{R}^{9}}\frac{\left\vert \triangle
u\right\vert ^{2}}{\left\vert \mathbf{x}_{1}-\mathbf{x}_{2}\right\vert
^{2}+\left\vert \mathbf{x}_{1}-\mathbf{x}_{3}\right\vert ^{2}+\left\vert 
\mathbf{x}_{2}-\mathbf{x}_{3}\right\vert ^{2}}dxdt \\
&\leqslant &C\int\limits_{0}^{2\pi }\int_{\mathbb{R}^{9}}\frac{\left\vert
\triangle u\right\vert ^{2}}{\left\vert \mathbf{x}_{1}-\mathbf{x}
_{2}\right\vert ^{2}}dxdt.
\end{eqnarray*}
due to the well-known change of variables $\mathbf{x}_{1}\rightarrow \frac{
\mathbf{x}_{1}-\mathbf{x}_{2}}{\sqrt{2}},\mathbf{x}_{2}\rightarrow \frac{
\mathbf{x}_{1}+\mathbf{x}_{2}}{\sqrt{2}}$ which is compatible with $
(-\triangle )$ and $(-\triangle +\left\vert x\right\vert ^{2}),$ we only
need to estimate
\begin{eqnarray}
&&\int\limits_{0}^{2\pi }\int_{\mathbb{R}^{9}}\frac{\left\vert \triangle
u\right\vert ^{2}}{\left\vert \mathbf{x}_{1}\right\vert ^{2}}dxdt
\label{example estimate} \\
&\leqslant &C\int\limits_{0}^{2\pi }\int_{\mathbb{R}^{9}}\frac{\left\vert
(-\triangle +\left\vert x\right\vert ^{2})u\right\vert ^{2}}{\left\vert 
\mathbf{x}_{1}\right\vert ^{2}}dxdt+C\int\limits_{0}^{2\pi }\int_{\mathbb{R}
^{9}}\frac{\left\vert \left\vert x\right\vert ^{2}u\right\vert ^{2}}{
\left\vert \mathbf{x}_{1}\right\vert ^{2}}dxdt  \notag \\
&\leqslant &C\left\Vert (-\triangle +\left\vert x\right\vert
^{2})f\right\Vert _{2}^{2}+C\int\limits_{0}^{2\pi }\int_{\mathbb{R}^{9}}
\frac{\left\vert \mathbf{x}_{1}\right\vert ^{4}\left\vert u\right\vert ^{2}}{
\left\vert \mathbf{x}_{1}\right\vert ^{2}}dxdt+C\int\limits_{0}^{2\pi
}\int_{\mathbb{R}^{9}}\frac{\left\vert \left\vert (\mathbf{x}_{2},\mathbf{x}
_{3})\right\vert ^{2}u\right\vert ^{2}}{\left\vert \mathbf{x}_{1}\right\vert
^{2}}dxdt  \notag \\
&\leqslant &C\left\Vert (-\triangle +\left\vert x\right\vert
^{2})f\right\Vert _{2}^{2}+C\left\Vert \nabla \hat{f}\right\Vert _{2}^{2}+E 
\notag \\
&\leqslant &C\left\Vert (-\triangle +\left\vert x\right\vert
^{2})f\right\Vert _{2}^{2}+E  \notag
\end{eqnarray}
where 
\begin{eqnarray*}
E &\leqslant &C\int\limits_{0}^{2\pi }\int_{\mathbb{R}^{9}}\frac{\left\vert
(-\triangle _{\mathbf{x}_{2},\mathbf{x}_{3}}+\left\vert (\mathbf{x}_{2},
\mathbf{x}_{3})\right\vert ^{2})u\right\vert ^{2}}{\left\vert \mathbf{x}
_{1}\right\vert ^{2}}dxdt+C\int\limits_{0}^{2\pi }\int_{\mathbb{R}^{9}}
\frac{\left\vert \triangle _{\mathbf{x}_{2},\mathbf{x}_{3}}u\right\vert ^{2}
}{\left\vert \mathbf{x}_{1}\right\vert ^{2}}dxdt \\
&\leqslant &C\int\limits_{0}^{2\pi }\int_{\mathbb{R}^{9}}\frac{\left\vert
(-\triangle _{\mathbf{x}_{2},\mathbf{x}_{3}}+\left\vert (\mathbf{x}_{2},
\mathbf{x}_{3})\right\vert ^{2})u\right\vert ^{2}}{\left\vert \mathbf{x}
_{1}\right\vert ^{2}}dxdt \\
&\leqslant &C\left\Vert (-\triangle +\left\vert x\right\vert
^{2})f\right\Vert _{2}^{2}
\end{eqnarray*}
due to lemma \ref{HermiteSoblev} and theorem \ref{main}.

Then it is easy to see that
\begin{eqnarray*}
B &=&36\int\limits_{0}^{2\pi }\int_{\mathbb{R}^{9}}a(x)\left\vert
u\right\vert ^{2}dxdt \\
&\leqslant &C\int\limits_{0}^{2\pi }\int_{\mathbb{R}^{9}}\frac{\left\vert
u\right\vert ^{2}}{\left\vert \mathbf{x}_{1}-\mathbf{x}_{2}\right\vert ^{2}}
dxdt \\
&\leqslant &C\left\Vert f\right\Vert _{2}^{2}
\end{eqnarray*}
because of theorem \ref{main} and change of variables.

The only term left over is 
\begin{equation*}
D=\int\limits_{0}^{2\pi }\int_{\mathbb{R}^{9}}a(x)V_{k}\func{Re}(u\overline{
u_{k}})dxdt.
\end{equation*}
A typical term in the sum reads
\begin{eqnarray*}
&&\left\vert \int\limits_{0}^{2\pi }\int_{\mathbb{R}^{9}}a(x)V_{1}\func{Re}
(u\overline{u_{1}})dx\right\vert \\
&\leqslant &\left( 4\int\limits_{0}^{2\pi }\int_{\mathbb{R}
^{9}}a(x)x_{1}^{2}\left\vert u\right\vert ^{2}dxdt\right) ^{\frac{1}{2}
}\left( \int\limits_{0}^{2\pi }\int_{\mathbb{R}^{9}}a(x)\left\vert \frac{
\partial }{\partial x_{1}}u\right\vert ^{2}dxdt\right) ^{\frac{1}{2}} \\
&\leqslant &C\left\Vert (-\triangle +\left\vert x\right\vert ^{2})^{\frac{1}{
2}}f\right\Vert _{2}^{2}
\end{eqnarray*}
using the same method as in estimate \ref{example estimate}.

Hence we conclude
\begin{equation*}
\int\limits_{0}^{2\pi }\int\limits_{\mathbb{R}^{3}}\left\vert u(\mathbf{x},
\mathbf{x,x,}t)\right\vert ^{2}d\mathbf{x}dt\leqslant C\left\Vert
(-\triangle +\left\vert x\right\vert ^{2})f\right\Vert _{2}^{2}
\end{equation*}
which is theorem \ref{collapsing variables}.

\begin{remark}
If we choose not to ignore $\int\limits_{0}^{2\pi }\int_{\mathbb{R}
^{9}}\triangle a\left\vert \nabla u\right\vert ^{2}dxdt$ and $
-\int\limits_{0}^{2\pi }\int_{\mathbb{R}^{9}}a(x)\left\vert \nabla
^{2}u\right\vert ^{2}dxdt$ in relation \ref{Morawetz}, then in fact we have
proven two additional Kato type smoothing estimates:
\begin{equation*}
\int\limits_{0}^{2\pi }\int_{\mathbb{R}^{9}}\frac{\left\vert \nabla
u\right\vert ^{2}}{(\left\vert \mathbf{x}_{1}-\mathbf{x}_{2}\right\vert
^{2}+\left\vert \mathbf{x}_{1}-\mathbf{x}_{3}\right\vert ^{2}+\left\vert 
\mathbf{x}_{2}-\mathbf{x}_{3}\right\vert ^{2})^{2}}dxdt\leqslant C\left\Vert
(-\triangle +\left\vert x\right\vert ^{2})f\right\Vert _{2}^{2}
\end{equation*}
and
\begin{equation*}
\int\limits_{0}^{2\pi }\int_{\mathbb{R}^{9}}\frac{\left\vert \nabla
^{2}u\right\vert ^{2}}{\left\vert \mathbf{x}_{1}-\mathbf{x}_{2}\right\vert
^{2}+\left\vert \mathbf{x}_{1}-\mathbf{x}_{3}\right\vert ^{2}+\left\vert 
\mathbf{x}_{2}-\mathbf{x}_{3}\right\vert ^{2}}dxdt\leqslant C\left\Vert
(-\triangle +\left\vert x\right\vert ^{2})f\right\Vert _{2}^{2}.
\end{equation*}
\end{remark}

\section{Appendix: Proof of Lemma \protect\ref{Antiderivative} / Computation
of the $L^{2}$ norms of the "antiderivatives" of Hermite functions}

\label{sec:odd and even computation}In this section, we prove lemma \ref
{Antiderivative} which yields the precise controlling constants. But we
shall first prove that there exits a $C>0$ s.t. 
\begin{equation*}
\left\Vert X_{k}\right\Vert _{2}^{2}\leqslant C,\text{ }\forall k
\end{equation*}
before we delve into the proof of lemma \ref{Antiderivative} which consists
of many special function techniques.

\subsection{Proof of the $L^{2}$ boundedness}

\begin{lemma}
\cite{Thangavelu}\label{creation and annihilation}We have the following
creation and annihilation relations
\begin{equation*}
(-\frac{d}{dx}+x)\tilde{h}_{k}(x)=\tilde{h}_{k+1}(x)
\end{equation*}
\begin{equation*}
(\frac{d}{dx}+x)\tilde{h}_{k}(x)=2k\tilde{h}_{k-1}(x)
\end{equation*}
where $\tilde{h}_{k}(x)=\frac{1}{c_{k}}h_{k}(x),$ and $c_{k}=(\frac{1}{
2^{k}k!\sqrt{\pi }})^{\frac{1}{2}}$ is the normalization constant, i.e. $
\tilde{h}_{k}(x)$ is the unnormalized Hermite function of degree $k$. In
this spirit, one has:
\begin{equation}
\tilde{h}_{k+1}(x)=-2\frac{d}{dx}\tilde{h}_{k}(x)+2k\tilde{h}_{k-1}(x)
\label{hermite recursive}
\end{equation}
or with the normalization factors
\begin{equation}
h_{k+1}(x)=-\sqrt{\frac{2}{k+1}}\frac{d}{dx}h_{k}(x)+\sqrt{\frac{k}{k+1}}
h_{k-1}(x).  \label{Recursive with normalization}
\end{equation}
\end{lemma}

We will only consider the even case
\begin{equation*}
V_{2k}=\frac{\left\Vert X_{2k}\right\Vert _{2}^{2}}{2}=\int\limits_{0}^{
\infty }(\int_{x}^{\infty }h_{2k}(t)dt)^{2}dx,
\end{equation*}
since the odd case is similar. Iterating relation \ref{Recursive with
normalization} yields 
\begin{equation}
h_{2k}(t)=\sum_{i=0}^{k-1}b_{i}\frac{d}{dt}h_{2k-1-2i}(t)+dh_{0}(t)
\label{even iteration}
\end{equation}
because 
\begin{eqnarray*}
h_{2k}(x) &=&-\sqrt{\frac{2}{2k}}\frac{d}{dx}h_{2k-1}(x)+\sqrt{\frac{2k-1}{2k
}}h_{2k-2}(x) \\
h_{2k-2}(x) &=&-\sqrt{\frac{2}{2k-2}}\frac{d}{dx}h_{2k-3}(x)+\sqrt{\frac{2k-3
}{2k-2}}h_{2k-4}(x) \\
&&... \\
h_{4}(x) &=&-\sqrt{\frac{2}{4}}\frac{d}{dx}h_{3}(x)+\sqrt{\frac{3}{4}}
h_{2}(x) \\
h_{2}(x) &=&-\sqrt{\frac{2}{2}}\frac{d}{dx}h_{1}(x)+\sqrt{\frac{1}{2}}
h_{0}(x).
\end{eqnarray*}

Therefore
\begin{eqnarray*}
\int_{0}^{\infty }\left( \int_{x}^{\infty }h_{2k}(t)dt\right) ^{2}dx
&\leqslant &2\int_{0}^{\infty }\left(
\sum_{i=0}^{k-1}b_{i}h_{2k-1-2i}(x)\right) ^{2}dx+2d^{2}\int_{0}^{\infty
}\left( \int_{x}^{\infty }h_{0}(t)dt\right) ^{2}dx \\
&\leqslant &2\int_{-\infty }^{\infty }\left(
\sum_{i=0}^{k-1}b_{i}h_{2k-1-2i}(x)\right) ^{2}dx+2d^{2}\int_{0}^{\infty
}\left( \int_{x}^{\infty }h_{0}(t)dt\right) ^{2}dx \\
&=&2(\sum_{i=0}^{k-1}\left\vert b_{i}\right\vert ^{2}+d^{2}\int_{0}^{\infty
}\left( \int_{x}^{\infty }h_{0}(t)dt\right) ^{2}dx)\text{,}
\end{eqnarray*}
where
\begin{eqnarray*}
\sum_{i=0}^{k-1}\left\vert b_{i}\right\vert ^{2} &=&\frac{2}{2k}+\frac{2k-1}{
2k}\frac{2}{2k-2}+\text{ }...\text{ }+\frac{2k-1}{2k}...\frac{2}{2} \\
&=&\frac{1}{k}+\frac{1}{k}\frac{2k-1}{2k-2}+\text{ }...\text{ }+\frac{1}{k}
\frac{2k-1}{2k-2}...\frac{3}{2} \\
&=&\frac{1}{k}\sum_{i=0}^{k-1}\left( \dprod\limits_{l=0}^{i}(1+\frac{1}{2k-2l
})\right) .
\end{eqnarray*}
Notice that
\begin{eqnarray*}
&&\ln \dprod\limits_{l=0}^{i}(1+\frac{1}{2k-2l}) \\
&\thicksim &\sum_{l=0}^{i}\frac{1}{2k-2l} \\
&\thicksim &\frac{1}{2}\ln \frac{k-1}{k-i}
\end{eqnarray*}
which implies
\begin{eqnarray*}
&&\sum_{i=0}^{k-1}\left\vert b_{i}\right\vert ^{2} \\
&\thicksim &\frac{1}{k}\sum_{i=0}^{k-1}(\frac{k-1}{k-i})^{\frac{1}{2}} \\
&\leqslant &\frac{1}{k}\sum_{i=0}^{k-1}(\frac{k}{k-i})^{\frac{1}{2}} \\
&\leqslant &C\frac{1}{k^{\frac{1}{2}}}\int_{0}^{k}(\frac{1}{k-x})^{\frac{1}{2
}}dx \\
&\leqslant &C
\end{eqnarray*}
i.e
\begin{equation*}
\left\Vert X_{2k}\right\Vert _{2}^{2}\leqslant C.
\end{equation*}

\begin{remark}
For the odd case, formula \ref{even iteration} will read
\begin{equation*}
h_{2k+1}(t)=\sum_{i=0}^{k-1}b_{i}\frac{d}{dt}h_{2k-2i}(t)+dh_{1}(t).
\end{equation*}
\end{remark}

\subsection{\protect\bigskip Proof of equalities \protect\ref{odd
computation} and \protect\ref{even computation}}

Below we will refer to the following lemmas as well as lemma \ref{creation
and annihilation}.

\begin{lemma}
\label{junk killing}Write the degree $k$ Hermite polynomial $e^{\frac{x^{2}}{
2}}\tilde{h}_{k}(x)$ as $H_{k}$, 
\begin{equation*}
H_{k}(x)=\sum_{i=0}^{[\frac{k}{2}]}\frac{k!}{i!(k-2i)!}(-1)^{i}(2x)^{k-2i}
\end{equation*}
then every polynomial $p(x)$ of degree $\leqslant i$ is a finite linear
combination of $H_{k},k\leqslant i,$
\begin{equation*}
p(x)=\sum_{k=0}^{i}(\int_{\mathbb{R}}\tilde{h}_{k}(x)p(x)e^{-\frac{1}{2}
x^{2}}dx)H_{k}(x)
\end{equation*}
In particular, given any polynomial $p(x)$ of degree$<k$, we have:
\begin{equation*}
\int_{\mathbb{R}}\tilde{h}_{k}(x)p(x)e^{-\frac{1}{2}x^{2}}dx=0
\end{equation*}
\end{lemma}

\begin{proof}
The first part of the statment is a well-known fact. To prove the second
part, one only needs to notice that $p(x)e^{-\frac{1}{2}x^{2}}$ is a finite
linear combination of $\tilde{h}_{i}(x)=H_{i}e^{-\frac{1}{2}x^{2}},i<k,$ and
then apply orthogonality.
\end{proof}

\begin{lemma}
\label{Laguerre}\bigskip \cite{Thangavelu}If we define the degree $k$
Lagurre polynomial of type $\alpha $ by 
\begin{equation*}
e^{-x}x^{\alpha }L_{k}^{\alpha }(x)=\frac{1}{k!}\frac{d^{k}}{dx^{k}}
(e^{-x}x^{\alpha }),
\end{equation*}
then
\begin{equation}
H_{2k+1}=(-1)^{k}2^{2k+1}k!L_{k}^{\frac{1}{2}}(x^{2})x.
\label{Thangavelu error}
\end{equation}
\end{lemma}

\begin{remark}
Formula \ref{Thangavelu error} is (1.1.53) in Thangavelu \cite{Thangavelu}.
He missed a factor $2$ on the right hand side. One can refer to page 1001 of 
\cite{Gradshteyn}.
\end{remark}

At this point we can give the proof of formula \ref{odd computation}

\subsubsection{Proof of the odd formula \protect\ref{odd computation}}

By relation \ref{hermite recursive}, we have
\begin{equation*}
\tilde{h}_{2k+1}(x)=-2\frac{d}{dx}\tilde{h}_{2k}(x)+4k\tilde{h}_{2k-1}(x)
\end{equation*}
and hence
\begin{equation*}
(\int_{-\infty }^{x}\tilde{h}_{2k+1}(t)dt)^{2}=4(\tilde{h}_{2k}(x))^{2}-16k
\tilde{h}_{2k}(x)\cdot \int_{-\infty }^{x}\tilde{h}_{2k-1}(t)dt+16k^{2}(
\int_{-\infty }^{x}\tilde{h}_{2k-1}(t)dt)^{2}
\end{equation*}
or with the normalization factors
\begin{equation*}
(X_{2k+1}(x))^{2}=\frac{2}{2k+1}(h_{2k}(x))^{2}-Junk(x)+\frac{2k}{(2k+1)}
(X_{2k-1}(x))^{2}
\end{equation*}
where 
\begin{equation*}
Junk(x)=\frac{16k}{(c_{2k+1})^{2}}\tilde{h}_{2k}(x)\cdot \int_{-\infty }^{x}
\tilde{h}_{2k-1}(t)dt.
\end{equation*}
So
\begin{equation*}
I_{2k+1}=\frac{2}{2k+1}+\int_{\mathbb{R}}Junk(x)dx+\frac{2k}{2k+1}I_{2k-1.}
\end{equation*}
where
\begin{equation*}
I_{2k+1}=\left\Vert X_{2k+1}\right\Vert _{2}^{2}
\end{equation*}
which is our target.

But $\tilde{h}_{2k-1}(t)$ $=$ $e^{-\frac{1}{2}x^{2}}(\dsum
\limits_{1}^{k}b_{2i-1}x^{2i-1})$, so $\int_{-\infty }^{x}\tilde{h}
_{2k-1}(t)dt=e^{-\frac{1}{2}x^{2}}(\dsum\limits_{i=0}^{k-1}l_{2i}x^{2i}),$
which implies $\int_{\mathbb{R}}Junk(x)dx=0$ by lemma \ref{junk killing}$.$
Hence
\begin{equation}
I_{2k+1}=\frac{2}{2k+1}+\frac{2k}{2k+1}I_{2k-1.}  \label{odd recursive}
\end{equation}

The equalities that 
\begin{equation*}
I_{1}=\int\limits_{-\infty }^{\infty }(\int_{-\infty }^{x}h_{1}(t)dt)^{2}dx=
\frac{1}{2\sqrt{\pi }}\int\limits_{-\infty }^{\infty }(\int_{-\infty
}^{x}2xe^{-\frac{1}{2}x^{2}}dt)^{2}dx=2
\end{equation*}
and relation \ref{odd recursive} tell us
\begin{equation*}
I_{2k+1}=2.
\end{equation*}

\subsubsection{Proof of the even formula \protect\ref{even computation}}

Applying relation \ref{Recursive with normalization} again, we have
\begin{eqnarray*}
&&V_{2k+2} \\
&=&\frac{1}{2}\frac{1}{k+1}+\frac{\sqrt{2}\sqrt{2k+1}}{k+1}
\int\limits_{0}^{\infty }h_{2k+1}(x)(\int_{x}^{\infty }h_{2k}(t)dt)dx+\frac{
2k+1}{2k+2}V_{2k} \\
&=&\frac{1}{2}\frac{1}{k+1}-\frac{\sqrt{2}\sqrt{2k+1}}{k+1}
\int\limits_{0}^{\infty }(\frac{d}{dx}\int_{x}^{\infty
}h_{2k+1}(t)dt)(\int_{x}^{\infty }h_{2k}(t)dt)dx+\frac{2k+1}{2k+2}V_{2k}
\end{eqnarray*}
Just as the odd case, we are concerned with the middle term and would like
to have an explicit formula for it. Integrating by parts once, we have 
\begin{eqnarray*}
&&\int\limits_{0}^{\infty }(\frac{d}{dx}\int_{x}^{\infty
}h_{2k+1}(t)dt)(\int_{x}^{\infty }h_{2k}(t)dt)dx \\
&=&-(\int_{0}^{\infty }h_{2k+1}(t)dt)(\int_{0}^{\infty
}h_{2k}(t)dt)+\int\limits_{0}^{\infty }(\int_{x}^{\infty
}h_{2k+1}(t)dt)h_{2k}(x)dx.
\end{eqnarray*}
Recall that we already know 
\begin{eqnarray*}
\int_{0}^{\infty }h_{2k}(t)dt &=&\frac{1}{2}\frac{2^{\frac{1}{2}
}2^{2k}\Gamma (\frac{1}{2}(2k)+\frac{1}{2})}{\sqrt{2^{2k}(2k)!\sqrt{\pi }}}
\\
&=&\frac{2^{k}\Gamma (k+\frac{1}{2})}{\sqrt{(2k)!\sqrt{\pi }}} \\
&=&\frac{2^{-k+\frac{1}{2}}(\sqrt{\pi })^{\frac{1}{2}}\Gamma (2k)}{\Gamma (k)
\sqrt{(2k)!}}
\end{eqnarray*}
from Feldheim \cite{Feldheim}, Busbridge \cite{BusBridge} and the well-known
formula for the gamma function
\begin{equation*}
\Gamma (z+\frac{1}{2})=\frac{2^{1-2z}\sqrt{\pi }\Gamma (2z)}{\Gamma (z)},
\end{equation*}
so we would like to compute $\int_{0}^{\infty }h_{2k+1}(t)dt$ and $
\int\limits_{0}^{\infty }(\int_{x}^{\infty }h_{2k+1}(t))h_{2k}(x)dx$. Using
lemma \ref{Laguerre}, we have 
\begin{eqnarray*}
\int_{0}^{\infty }h_{2k+1}(t)dt &=&\frac{1}{c_{2k+1}}(-1)^{k}2^{2k+1}k!
\int_{0}^{\infty }L_{k}^{\frac{1}{2}}(x^{2})e^{-\frac{x^{2}}{2}}xdx \\
&=&\frac{(-1)^{k}2^{2k+1}k!}{2c_{2k+1}}\int_{0}^{\infty }L_{k}^{\frac{1}{2}
}(u)e^{-\frac{u}{2}}du \\
&=&\frac{(-1)^{k}2^{2k+1}k!}{c_{2k+1}}\sum_{i=0}^{k}\binom{\frac{1}{2}+i-1}{i
}(-1)^{k-i} \\
&=&\frac{2^{k+1}k!}{\sqrt{2(2k+1)!\sqrt{\pi }}}\sum_{i=0}^{k}\binom{\frac{1}{
2}+i-1}{i}(-1)^{i}
\end{eqnarray*}
where the integral part, which has been worked out in page 809 of \cite
{Gradshteyn}, is that
\begin{equation*}
\int_{0}^{\infty }L_{k}^{\alpha }(u)e^{-\beta u}du=\sum_{i=0}^{k}\binom{
\alpha +i-1}{i}\frac{(\beta -1)^{k-i}}{\beta ^{k-i+1}}.
\end{equation*}
Hence
\begin{eqnarray*}
&&(\int_{0}^{\infty }h_{2k+1}(t)dt)(\int_{0}^{\infty }h_{2k}(t)dt) \\
&=&\frac{2^{-k+\frac{1}{2}}(\sqrt{\pi })^{\frac{1}{2}}\Gamma (2k)}{\Gamma (k)
\sqrt{(2k)!}}\frac{2^{k+1}k!}{\sqrt{2(2k+1)!\sqrt{\pi }}}\sum_{i=0}^{k}
\binom{\frac{1}{2}+i-1}{i}(-1)^{i} \\
&=&\frac{1}{\sqrt{2k+1}}\sum_{i=0}^{k}\binom{\frac{1}{2}+i-1}{i}(-1)^{i} \\
&=&\frac{1}{\sqrt{2k+1}}\sum_{i=0}^{k}\binom{-\frac{1}{2}}{i},
\end{eqnarray*}
due to the identity 
\begin{equation*}
\binom{\alpha }{k}=\binom{k-\alpha -1}{k}(-1)^{k}.
\end{equation*}

For the last term, we have 
\begin{eqnarray*}
\int\limits_{0}^{\infty }(\int_{x}^{\infty }h_{2k+1}(t))h_{2k}(x)dx &=&
\frac{1}{2}\int\limits_{-\infty }^{\infty }(\int_{x}^{\infty
}h_{2k+1}(t)dt)h_{2k}(x)dx \\
&=&\frac{1}{2}\frac{(\text{the coefficient of }x^{2k}e^{-\frac{x^{2}}{2}
\text{ }}\text{in }\int_{x}^{\infty }h_{2k+1}(t)dt)}{(\text{the coefficient
of }x^{2k}e^{-\frac{x^{2}}{2}}\text{ in }h_{2k}(x))} \\
&=&\frac{1}{2}\frac{\frac{2^{2k+1}}{c_{2k+1}}}{\frac{2^{2k}}{c_{2k}}} \\
&=&\frac{1}{\sqrt{2}\sqrt{2k+1}},
\end{eqnarray*}
via lemma \ref{junk killing}.

At long last we have 
\begin{equation*}
V_{2k+2}=-\frac{1}{2k+2}+\frac{\sqrt{2}}{k+1}\sum_{i=0}^{k}\binom{-\frac{1}{2
}}{i}+\frac{2k+1}{2k+2}V_{2k}.
\end{equation*}
Since
\begin{equation*}
\frac{1}{k+1}\sum_{i=0}^{k}\binom{-\frac{1}{2}}{i}+\frac{2k+1}{2k+2}
\sum_{i=0}^{k}\binom{\frac{1}{2}}{i}=\sum_{i=0}^{k+1}\binom{\frac{1}{2}}{i},
\end{equation*}
a straight forward induction gives us formula \ref{even computation}. This
concludes the proof of lemma \ref{Antiderivative}.

\end{document}